\documentclass[10 pt, conference]{ieeeconf}
\usepackage[utf8]{inputenc}

\usepackage[bookmarks=true]{hyperref}
\usepackage{xcolor}
\usepackage{graphicx}
\usepackage{dsfont} 
\usepackage{amsmath}
\usepackage{amssymb}
\usepackage{mathtools} 
\newtheorem{remark}{Remark}

\usepackage{url}

\usepackage{color}
\usepackage{subcaption}
\usepackage[utf8]{inputenc}
\usepackage[T1]{fontenc}
\usepackage[T1]{fontenc}
\IEEEoverridecommandlockouts
\usepackage{cite}

\usepackage{cleveref}

\newtheorem{assumption}{Assumption}
\newtheorem{definition}{Definition}
\newtheorem{lemma}{Lemma}
\newtheorem{theorem}{Theorem}

\usepackage{comment}
\usepackage{mathrsfs}
\usepackage[textfont=md,font=footnotesize]{caption}
\DeclareMathOperator*{\argmin}{arg\,min}

\newcommand{\doi}[1]{\href{http://dx.doi.org/#1}{\normalsize{\textsc{doi:}}~\nolinkurl{#1}}}
\newcommand{\arxiv}[1]{\href{http://arxiv.org/abs/#1}{\normalsize{\textsc{arxiv:}}~\nolinkurl{#1}}}
\newcommand{\todo}[1]{{\normalsize{\textbf{({\color{red}TODO:}\ #1)}}}}

\renewcommand{\epsilon}{\varepsilon}
\renewcommand{\phi}{\varphi}

\newcommand{\R}{\mathbb{R}}

\newcommand{\N}{\mathbb{N}}

\let\originalleft\left
\let\originalright\right
\renewcommand{\left}{\mathopen{}\mathclose\bgroup\originalleft}
\renewcommand{\right}{\aftergroup\egroup\originalright}

\def\clap#1{\hbox to 0pt{\hss#1\hss}}

\newcommand{\norm}[1]{\left\lVert #1\right\rVert}

\newcommand{\set}[1]{\left\{ #1\right\}}

\let\sp=\sqparen

\DeclareMathAlphabet{\mathpzc}{OT1}{pzc}{m}{it}

\newcommand{\fcnote}[1]%
    {\textcolor{orange}{\textbf{FC: #1}}}
\newcommand{\twnote}[1]%
    {\textcolor{cyan}{\textbf{TW: #1}}}
\newcommand{\aanote}[1]%
    {\textcolor{blue}{\textbf{AA: #1}}}

\title{Learning Min-norm Stabilizing Control Laws for \\
Systems with Unknown Dynamics}

\begin{document}

\author{Tyler Westenbroek$^{*1}$\thanks{$^*$ Indicates equal contribution.}, Fernando Casta\~neda$^{*2}$, Ayush Agrawal$^{2}$, S. Shankar Sastry$^{1}$, Koushil Sreenath$^{2}$
\thanks{$^{1}$Department of Electrical Engineering and Computer Sciences, University of California at Berkeley, USA.}%
\thanks{$^2$Department of Mechanical Engineering, University of California at Berkeley, USA.}
\thanks{The work of Fernando Casta\~neda received the support of a fellowship (code LCF/BQ/AA17/11610009) from ”la Caixa” Foundation (ID 100010434). This work was partially supported through National Science Foundation Grants CMMI-1931853 and CMMI-1944722 and by HICON-LEARN (design of HIgh CONfidence LEARNing-enabled systems), Defense Advanced Research Projects Agency award number FA8750-18-C-0101, and Provable High Confidence Human Robot Interactions, Office of Naval Research award number N00014-19-1-2066.}
}

\maketitle
\begin{abstract}
This paper introduces a framework for learning a minimum-norm stabilizing controller for a system with unknown dynamics using model-free policy optimization methods. The approach begins by first designing a Control Lyapunov Function (CLF) for a (possibly inaccurate) dynamics model for the system, along with a function which specifies a minimum acceptable rate of energy dissipation for the CLF at different points in the state-space. Treating the energy dissipation condition as a constraint on the desired closed-loop behavior of the real-world system, we use penalty methods to formulate an unconstrained optimization problem over the parameters of a learned controller, which can be solved using model-free policy optimization algorithms using data collected from the plant. We discuss when the optimization learns a stabilizing controller for the real world system and derive conditions on the structure of the learned controller which ensure that the optimization is strongly convex, meaning the globally optimal solution can be found reliably. We validate the approach in simulation, first for a double pendulum, and then generalize the framework to learn stable walking controllers for underactuated bipedal robots using the Hybrid Zero Dynamics framework. By encoding a large amount of structure into the learning problem, we are able to learn stabilizing controllers for both systems with only minutes or even seconds of training data.
\end{abstract}

\section{Introduction}
\label{sec:intro}
Recently, the literature has displayed a renewed interest in data-driven methods for controller design \cite{taylor2019episodic, berkenkamp2017safe, marco2016automatic, akametalu2014reachability}. Much of this excitement has been driven by recent advances in the model-free reinforcement learning literature \cite{hwangbo,levine_end_to_end}. Despite their generality, model-free policy optimization methods are known to suffer from poor sample complexity, as they generally are unable to take advantage of known structure in the control system. This paper bridges the gap between model-based and model-free design paradigms by embedding Lyapunov-based design techniques into a model-free reinforcement learning problem. By encoding basic information about the structure of the system into the learning problem through a Control Lyapunoc Function (CLF), our approach is able to learn optimal stabilizing controllers for highly uncertain systems with as little as seconds or a few minutes of data. 

Specifically, the paper proposes a framework for learning a min-norm stabilizing control law for an unknown system using model-free policy optimization techniques. Our approach begins by first designing a CLF for a nominal dynamics model of the system alongside a function which specifies the desired rate of convergence for the closed-loop system. To impose this desired behavior on the real world control system, we then formulate a continuous-time optimization problem over the parameters of a learned controller which treats the energy dissipation condition as a constraint. The cost function for the optimization encourages choices of parameters which minimize control effort, but uses a penalty term to ensure that the dissipation constraint is satisfied, if possible, when the penalty term is chosen to be large enough. The terms in the optimization depend on the dynamics of the unknown system, but discrete-time approximations to the problem can be solved using policy-optimization algorithms and data collected from the plant. In general, the problem may be non-convex, but when the learned controller is linear in its parameters the problem becomes (strongly) convex, meaning the globally optimal solutions for the problem can be found using standard policy gradient \cite{sutton2018reinforcement} or random search techniques \cite{mania2018simple}. 

To demonstrate the utility of the proposed framework, we apply the method in simulation to a double pendulum and a high-dimensional model of a bipedal robot. For the double pendulum example, the learned controller is comprised of a linear combination of radial basis functions so that the convexity result discussed above applies, and we demonstrate empirically that the learned controller is able to closely match the true min-norm controller performance. The walking example demonstrates how to extend our results in the body of the paper to encompass the Hybrid Zero Dynamics framework as in \cite{ames2014rapidly}. For this high-dimensional system, a feed-forward neural network is used for the learned controller. While we cannot guarantee that the optimal set of parameters is found, the learned controller still produces a stable walking motion in the face of high model uncertainty.

\subsection{Related Work}
CLF-based controllers \cite{ARTSTEIN19831163,SONTAG1989117} have been proved to be effective for a wide variety of complex robotic tasks, such as bipedal walking \cite{galloway_clf},\cite{ames2014rapidly}, manipulation \cite{ames_clf_manipulation} and multi-agent coordination \cite{Ogren_clf}. In \cite{galloway_clf} and \cite{ames_clf_manipulation} quadratic programs (CLF-QP), which integrate the CLF condition as a constraint, are used to get optimal min-norm stabilizing controllers. The CLF-QP is solved online and additional constraints, such as input saturation, can be added.

However, the dynamics of many real-world systems have nonlinearities that might be difficult to model correctly and/or physical parameters which could be difficult to identify. Input-to-state stability has been used to tackle this problem in \cite{HESPANHA20082735},\cite{SONTAG199524}. Also, adaptive\cite{ACC2015_L1_CLF} and robust\cite{robust_clf_battilotti},\cite{RSS2015_RobustCLF} versions of CLF-based controllers have been developed in recent years. However, these approaches sometimes fail to account for the correct amount of uncertainty due to the typical assumptions they make on the uncertainties' structures and bounds.

Our work most closely aligns with recent research that use data-driven approaches to tackle the issue of model uncertainty in nonlinear controllers. Our works builds on \cite{westenbroek2019feedback, castaneda20a}, where reinforcement learning is used to account for uncertainty when performing feedback linearization of nonlinear systems. In contrast to recently proposed approaches \cite{taylor2019episodic, choi2020reinforcement} which focus on learning the uncertain terms in a CLF-QP in order to indirectly improve the optimization-based controller, the framework proposed in this paper directly learns the optimal stabilizing controller. By directly learning the desired controller, our approach removes the need for solving a real-time optimization problem involving a potentially complex learned component, which may take a non-trivial amount of time to process during real-time applications. On hardware, CLF-based controllers frequently need to be updated at frequencies exceeding 1000 Hertz to maintain the stability of the system \cite{ames2014rapidly}, placing strict timing requirements on the rate at which the CLF-based controller must be updated. Thus, we hypothesize that the direct approach may have meaningful advantages in applications, and future work will seek to validate this claim.

\subsection{Organization}
The rest of the paper is organized as follows. Section \ref{sec:clfs} revisits Control Lyapunov Functions. Section \ref{sec:learn-min-norm} presents the proposed learning problem, develops our theoretical guarantees, and then demonstrates how the discrete-time approximations to the problem can be solved using reinforcment learning in an approach similar to \cite{westenbroek2019feedback,castaneda20a}. In Section \ref{sec:examples} the proposed method is used to stabilize a double pendulum and the walking gait of an underactuated nonlinear bipedal robot. Finally, Section \ref{sec:conclusion} provides concluding remarks.
\subsection{Notation and Terminology}
We let $\R_{\geq 0} = \set{x \in \mathbb{R} \colon x \geq 0}$ denote the closed right half plane. We say that a function $V \colon \mathbb{R}^n \to \mathbb{R}_{\geq 0}$ is \emph{positive definite} if $V(0) = 0$ and $V(x)>0$ if $x \neq 0$. We further say that $V$ is \emph{radially unbounded} if $V(x) \to \infty$ as $\| x\| \to \infty$. If it exists, the gradient of a function $V \colon \mathbb{R}^n \to \mathbb{R}$ at the point $x \in \mathbb{R}^n$ is denoted by the row vector $\nabla V(x) \in \mathbb{R}^{1 \times n}$.

\section{Control Lyapunov Functions}
\label{sec:clfs}
Throughout the paper we will consider nonlinear control-affine systems of the form
\begin{equation}\label{eq:vfield}
    \dot{x} = f(x) + g(x)u,
\end{equation}
where $x \in \R^n$ is the state and $u \in \mathbb{R}^m$ the input. The mappings $f \colon \R^n \to \R^n$ and $g \colon \R^n \to \R^{n \times m}$ are assumed to be locally Lipschtiz continuous with $f(0)=0$.


\begin{definition}\vspace{0.3cm}
We say that a continuously differentiable, positive definite, radially unbounded function $V \colon \R^n \to \R_{\geq0}$  is a \emph{Control Lyapunov Function} (CLF) for \eqref{eq:vfield} with positive definite energy dissipation rate $\sigma \colon \mathbb{R}^n \to \mathbb{R}_{\geq 0}$ if for each $x \in \mathbb{R}^n \setminus \set{0}$
\begin{equation}\label{eq:clf_def}
    \inf_{u \in \mathbb{R}^m} \nabla V(x) \cdot [f(x) + g(x)u] \leq -\sigma(x).
\end{equation} \vspace{0.05cm}
\end{definition}

 It is well-known that if the above conditions are satisfied then the system is asymptotically controllable in the sense that the state can be driven to the origin asymptotically for every initial condition \cite{SONTAG1989117}. For many physical systems it is desirable to find a locally Lipschitz continuous feedback rule $u \colon \R^n \to \R^m$ so that for each $x \in \R^n \setminus \set{0}$
\begin{equation}\label{eq:closed_loop_stable}
    \nabla V(x) \cdot [f(x) + g(x)u(x)] \leq -\sigma(x)
\end{equation}
and the closed loop system is asymptotically stable. It should be noted that not all systems which satisfy \eqref{eq:clf_def} admit such a controller \cite{sastry2013nonlinear}, but a number of important systems such as the ones considered in this document are continuously stabilizable.
One popular choice of control law which satisfies the dissipation constraint \eqref{eq:closed_loop_stable} is the min-norm control law $u^* \colon \R^n \to \R^m$ which is defined point-wise by: 
\begin{align} \label{eq:min_norm_controller} 
    u^*(x) &= \argmin_{u \in \R^m} \|u\|_2^2 \nonumber \\
         & \text{s.t.} \     \nabla V(x) \cdot [f(x) + g(x)u] \leq -\sigma(x)
\end{align}
At every point, this controller selects the smallest input which ensures that the CLF decays at the desired rate. If $V$ is a CLF for the system, a sufficient condition for $u^*$ to be locally Lipschitz continuous is that $f$, $g$ and the gradient of $V$ are each locally Lipschitz continuous \cite{freeman2008robust}. Moreover, when there are no constraints on the input, one can derive a closed-form expression for the controller (see e.g. \cite{ames2014rapidly}).  
Moreover, letting 
\begin{equation}
    a(x)= \nabla V(x) \cdot f(x) + \sigma(x) \text{ and } b(x) = \nabla V(x) \cdot g(x),
\end{equation}
the min-norm has the following closed-form representation:
\begin{equation}
    u^*(x) = \begin{cases}
    - \frac{a(x)b(x)}{b(x)^{T}b(x)} & \text{ if } a(x) >0 \\
    0 & \text{ if } a(x) \leq 0
    \end{cases}
\end{equation}
However, one advantage of formulating the min-norm controller as a point-wise optimization as in \eqref{eq:min_norm_controller} is that the optimization can easily incorporate bounds on the allowable control efforts for the system by restricting $u \in U \subset \R^{m}$. This is important in many applications where the actuators of the system have physical limitations. Our theoretical results rely on the input being unconstrained, however, in practice such constraints can be added by restricting the range-space of the learned controller.

\section{Learning Min-norm Stabilizing Controllers}
\label{sec:learn-min-norm}
\subsection{Learning a Min-norm Stabilizing Controller for a System with Unknown Dynamics} \label{subsec:learn-min-norm}

Despite the wide-spread utility of the CLF-based controllers introduced in the previous section, the primary drawback of these methods is that they require an accurate dynamics model to implement. Our present objective is to learn a min-norm stabilizing controller for the plant with unknown dynamics
\begin{equation}\label{eq:plant}
    \dot{x} = f_p(x) + g_p(x)u,
\end{equation}
while ensuring that the learned controller adheres to the dissipation constraint imposed by some candidate CLF $V \colon \R^n \to \R_{\geq 0}$ and associated decay rate $\sigma \colon \R^n \to \R_{\geq 0}$.  

We will focus on learning the min-norm controller for the plant on a compact subset of the state-space. Specifically, we will focus on learning the min-norm stabilizing controller for the system on the set
\begin{equation}
    W^c \coloneqq \set{x \in \R^n \colon V(x) \leq c},
\end{equation}
where $c>0$ is a design parameter. 

We will make the following technical assumptions throughout the paper unless otherwise specified:
\begin{assumption}\vspace{0.3cm}\label{ass:lispschitz_stuff}
The components $f_p$, $g_p$, $\sigma$ and $\nabla V$ are each locally Lipschitz continuous.\vspace{0.3cm}
\end{assumption}
\begin{assumption}\label{ass:clf_existence}
There exists a locally Lipschitz continuous control law $\tilde{u}_p \colon \R^n \to \R^m$ such that for each $x \in W^{c}$
\begin{equation}\label{eq:plant_clf_constraint}
    \nabla V(x) [f_p(x) + g_p(x)\tilde{u}_p(x)] \leq -\sigma(x).
\end{equation}
\end{assumption}\vspace{0.3cm}
\begin{remark}
Assumption \ref{ass:clf_existence} ensures that $V$ is a true CLF for the system with associated dissipation rate $\sigma$. While our approach does not explicity require a nominal dynamics model for the plant, in practice, our candidate CLF for the plant is constructed using a nominal dynamics model
\begin{equation}\label{eq:nominal_dynamics}
    \dot{x} = f_m(x) + g_m(x)u,
\end{equation}
which incorporates any information we have about the plant, but may be inaccurate due to nonlinearities which are difficult to model or dynamics parameters which are challenging to identify. However, despite model mismatch between \eqref{eq:plant} and \eqref{eq:nominal_dynamics}, we can often design a CLF for the model and reasonably expect it to also be a CLF for the plant. For example, our two numerical examples systematically construct CLF's for the unknown system using feedback-linearizing coordinates. In these examples Assumption \ref{ass:clf_existence} is tantamount to knowing the relative degree of the system, a rather mild structural assumption. 
\end{remark}

Since we do not know the terms in \eqref{eq:plant}, we now propose a method to learn a stabilizing CLF-based controller for the system using data collected from the plant. Under the preceding assumptions, and recalling our discussion from Section \ref{sec:clfs},
we know that there is a well-defined control law $u_p^* \colon W^{c} \to \R_{\geq 0}$ which asymptotically stabilizes the plant on $W^{c}$ and is given point-wise by
\begin{align*}
    u_p^*(x) & = \argmin_{u \in \R^m} \norm{u}_2^2 \\ 
    & \text{s.t.}  \ \nabla V(x)[f_p(x) + g_p(x)u] \leq -\sigma(x).
    \vspace{0.3cm}
\end{align*}
We will denote our learned approximation for $u_p^*$ by $\hat{u} \colon \R^n \times \Theta \to \R^m$. For each choice of parameter $\theta \in \Theta \subset \R^K$ the control law $\hat{u}(\cdot, \theta) \colon \R^n \to \R^m$ defines the learned control law supplied to the plant, with $\Theta \subset \mathbb{R}^K$ a convex set of allowable learned parameters. It is assumed that $\hat{u}$ is locally Lipschitz continuous in its first argument and continuously differentiable in its second argument. Common function approximators such as feed-forward neural networks, radial basis functions or bases of polynomials can be used to construct the learned controller.

\begin{remark}
In general, the learned controller can incorporate information from a nominal dynamics model by giving it the structure 
\begin{equation}\label{eq:augmented_controller}
    \hat{u}(x,\theta) = u_m(x) + \delta u(x,\theta),
\end{equation}
where $u_m$ is a nominal model-based controller and $\delta u \colon \R^n \times \R^K \to \R^m$ is the learned component. 
\end{remark}

Next, in order to find parameters for the learned controller which satisfy the dissipation constraint \eqref{eq:plant_clf_constraint}, we will solve optimizations over the parameters of the learned controller of the form
\begin{equation}\label{eq:clf_opt}
(\textbf{P}_\lambda): \min_{\theta \in \Theta} L_\lambda(\theta),
\end{equation}
where for each $\lambda \in \R_{\geq0}$ we define the loss function
\begin{equation}
 L_\lambda(\theta) =  E_{x \sim X} \left[\| u(x,\theta)\|_2^2 + \lambda H(\Delta(x,\theta))\right],
\end{equation}
 where $X$ is the uniform probability distribution over $W^{c}$, the mapping $\Delta \colon \R^n \times \Theta \to \R$ is defined by
\begin{equation}\label{eq:constraint_cost}
  \Delta(x,\theta) = \nabla V(x)[f_p(x) + g_p(x)\hat{u}(x,\theta)] + \sigma(x)
\end{equation}
and finally  $H \colon \R \to \R_{\geq 0}$ is defined for each $y \in \R$ by
\begin{equation}
    H(y) = \begin{cases}
     y & \text{ if } y \geq 0, \\
    0 & \text{ if } y < 0.
    \end{cases}
\end{equation}
The first term in the loss $L_\lambda$ encourages small control efforts while the second term penalizes violations of the CLF dissipation constraint, with $\lambda \in \R_{\geq 0}$ used to control the magnitude of the penalty.  While we do not know $\Delta(x,\theta)$ \emph{a priori}, we can measure this quantity by applying the control $\hat{u}(x,\theta)$ to the plant at the point $x$ and measuring the resulting time derivative of $V$. Then, equation \eqref{eq:constraint_cost} can be used to compute the desired quantity. Thus, any stochastic optimization algorithm can be used to solve $\textbf{P}_\lambda$ by running experiments to evaluate the terms in $L_\lambda$. We will discuss this in further detail when we present practical approaches for solving $\textbf{P}_\lambda$ below. 
\begin{remark}
The uniformity of the distribution $X$ ensures that all points in $W^c$ are considered when optimizing over the parameters of $\hat{u}$. This requirement is seen to be analogous to the persistency of excitation conditions which are common in the adaptive control literature \cite{sastry1989adaptive}. In the proofs of the following theoretical results this condition ensures that each component of the learned controller is activated sufficiently during the learning process.  
\end{remark}

\subsection{Theoretical Results}\label{subsec:theory}
We now study how the solution set of $\textbf{P}_\lambda$ changes as the penalty parameter $\lambda$ is increased and derive conditions under which the problem is convex, meaning that it can be solved reliably to global optimality using iterative gradient-based optimization algorithms. To simplify the statement of our results, for each $\lambda \in \R_{\geq 0}$ we define
\begin{equation}
    S_\lambda \colon = \set{\theta \in \Theta \colon \theta \in \argmin_{\theta \in \Theta} L_\lambda(\theta)}
\end{equation}
to capture the set of global minimizers for $\textbf{P}_\lambda$. We also define
\begin{equation}\label{eq:safe_values}
   \Xi \colon = \set{\theta \in \Theta \colon \Delta(x,\theta) \leq 0 , \ \forall x\in W^c} \subset \Theta
\end{equation}
to be the set of parameters for which the corresponding learned controller satisfies the desired CLF dissipation constraint at every point in $W^c$. Next, we present our theoretical results in Lemma \ref{lemma:penalty} and Theorems \ref{thm:learn_min_norm} and \ref{thm:convex}, whose proofs can be found in the Appendix.

First, we compare the sets $\Xi$ and $S_\lambda$ as the penalty term $\lambda$ is increased:
\begin{lemma}\vspace{0.1cm} \label{lemma:penalty}
Assume that $\Xi$ is non-empty so that there exists at least one choice of learned parameters which satisfy the desired CLF constraint. Then there exists $\bar{\lambda} \in \R_{\geq0}$ such that for each $\lambda > \bar{\lambda}$ all global optimizers of $\textbf{P}_\lambda$ also satisfy the dissipation constraint, namely, $S_\lambda \subset \Xi$.
\vspace{0.1cm}
\end{lemma}

In other words, if the penalty parameter $\lambda \in \R_{\geq0}$ is chosen to be large enough then $\textbf{P}_\lambda$ recovers the set of learned parameters which stabilize the plant and satisfy the CLF constraint. Note that if $\theta^* \in \Xi$ is one such choice of parameters then it must be the case that $\mathbb{E}_{x \sim X} \left[\lambda H(\Delta(x,\theta^*))\right] =0 $. Thus, when $\Xi$ is non-empty and $\lambda$ is chosen to be large enough the minimizers of $\textbf{P}_\lambda$ are selected by the set of parameters which minimize the term $\mathbb{E}_{x\sim X}\left[\norm{u(x,\theta)}_2^2\right]$, which is the average control effort exerted over the state-space by the corresponding learned controller. By definition, the min-norm stabilizing controller $u_p^*$ minimizes the control effort needed to satisfy the CLF dissipation constraint at every point in the state-space. Thus, if $\lambda$ is large enough and $u_p^*$ is in the space of learned controllers spanned by $\hat{u}$, it must be recovered by the optimization:

\begin{theorem}\label{thm:learn_min_norm}\vspace{0.1cm}
Assume that there exists $\bar{\theta} \in \Theta$ such that $\hat{u}(x,\bar{\theta}) =u_p^*(x)$ for each $x \in W^c$. Then there exists $\bar{\lambda} \in \R_{\geq0}$ such that for each $\lambda>\bar{\lambda}$ and $\theta^* \in S_\lambda$ we have $\hat{u}(x,\theta^*) = u_p^*(x)$ for each $x \in W^c$. \vspace{0.1cm}
\end{theorem}

However, the family of optimization problems we have formulated over the parameters of the learned controller will generally be non-convex, meaning that we cannot efficiently find their globally optimal solutions. Thus, we seek conditions under which $\textbf{P}_\lambda$ becomes convex so that we can reliably find its global minimizers using iterative methods. Towards this end we will now assume that 
\begin{equation}\label{eq:lin_param}
    \hat{u}(x,\theta) = \sum_{k = 1}^{K} \theta_k u_k(x),
\end{equation}
where $\set{u_k}$ is a set of locally Lipschitz continuous mappings from $\R^n$ to $\R^m$ and $\theta_k$ is the $k$-th entry of the learned parameter\footnote{Alternatively, one could also assume that the learned controller is of the form $\hat{u} = u_m(x) + \sum_{k = 1}^{K} \theta_k u_k(x)$ if the system designer wishes to augment a known model-based controller as in \eqref{eq:augmented_controller}. The statement and proof of Theorem \ref{thm:convex} go through with minor modifications in this case.}. Linearity assumptions of this sort are common in convergence proofs found in both the adaptive control and reinforcement learning literature. In the statement of the following result we informally view each basis element $u_k$ as a subset of $C(W^c, \R^m)$, the vector space of continuous functions from $W^c$ to $\R^m$. 

\begin{lemma}\vspace{0.1cm}\label{lemma:convex}
Assume that $\hat{u}$ is of the form \eqref{eq:lin_param} and that the set $\set{u_k}_{k=1}^K$ is linearly independent on $C(W^{c},\R^m)$. Then for each $\lambda \in \R_{\geq 0}$ the loss $L_\lambda$ is strongly convex. \vspace{0.1cm}
\end{lemma}

This leads to the main theoretical result of this paper, which follows from an immediate application of Theorem \ref{thm:learn_min_norm} and Lemma \ref{lemma:convex}.

\begin{theorem}\label{thm:convex} \vspace{0.1cm}
Assume that $u$ is of the form \eqref{eq:lin_param} and that the set $\set{u_k}_{k=1}^K$ is linearly independent on $C(W^{c},\R^m)$. Further assume that $\Theta \subset \R^K$ is convex and that there exists $\theta^* \in \Theta$ such that $\hat{u}(x,\theta^*) = u^*_p(x)$ for each $x \in W^{c}$. Then there exists $\bar{\lambda} \in \R_\lambda \geq 0$ such that for each $\lambda >\bar{\lambda}$ the problem $\textbf{P}_\lambda$ is a strongly convex optimization problem with $\theta^*$ its unique global minimizer. \vspace{0.1cm}
\end{theorem}

In practice, since we do not have a parametric model for the system, it is unlikely that there exists a set of learned parameters which exactly reconstructs the true min-norm controlller for the plant. However, many model-free learning schemes \cite{lewis2009reinforcement} make use of function approximation schemes which can recover a continuous function up to a desired level of accuracy, if enough terms are included in the basis of features. It is a matter for future work to show that we can leverage these results to learn $u^*(x)$ up to a pre-specified degree of accuracy.  However, in practice the number of elements required in such an expansion can quickly become prohibitively large as the dimension of the state grows. Thus, for high dimensional systems, such as the bipedal robots we consider below, practical implementations may require the use of more compactly represented function approximation schemes, such as multi-layer feed-forward neural networks, which can also approximate continuous functions to a desired degree of accuracy (Universal Approximation Theorem \cite{Cybenko1989, HORNIK_approximation}), but lead to non-convexities in our optimization problem \eqref{eq:lin_param}. 

\subsection{Solving Discrete-time Approximations with Reinforcement Learning}\label{subsec:rl}
 Many real-world systems have digital sensors and actuators which can only be updated at some maximum frequency, meaning we can only obtain finite difference approximations of $\dot{V}$ when different control signals are applied to the plant. Thus, in this section we introduce discrete-time approximations to $\textbf{P}_\lambda$ and discuss how they can be solved with standard reinforcement learning algorithms \cite{schulman2017proximal, levine_sac}. Our description of this process will be brief, since the approach is similar to the one described in \cite{westenbroek2019feedback}.

For the reinforcement learning problem we will assume that the control supplied to the plant can only be updated at a fixed minimum sampling period $\Delta t>0$. We will let $t_k = k \times \Delta t$ for each $k \in \mathbb{N}$ denote the set of sampling intervals. When the control $\hat{u}(x,\theta) \in \R^m$ is applied over the interval $[t_k, t_{k+1}]$ a Taylor expansion can be used to show that 
\begin{equation}
    \Delta(x,\theta) =  \underbrace{\frac{V(x(t_{k+1})) - V(x(t_{k})) }{\Delta t} + \sigma(x(t_k))}_{\colon = \tilde{\Delta} (x,\theta)} + O(\Delta t^2). 
\end{equation}
Thus for small $\Delta t$ we use the loss
\begin{equation}
    \tilde{l}_\lambda(x,\theta) =  \| u(x,\theta)\|_2^2 + \lambda H(\tilde{\Delta}(x,\theta)).
\end{equation}
We use this approximate pointwise loss to define the following reinforcement learning problem, which serves as an approximation to $\textbf{P}_\lambda$:
\begin{align}\label{eq:rl_prob}
    \tilde{\textbf{P}}_\lambda \colon  \min_{\theta \in \Theta} &\ E_{x_0 \sim X} \left[\sum_{k = 0}^N \tilde{l}(x_k,\theta)\right]\\ \nonumber
    \text{s.t. } x_{k+1} &= x_k + \int_{t_k}^{t_{k+1}}\left[f(x(t)) +g(x(t))u_k\right]dt\\ \nonumber
    u_k &= \hat{u}(x_k, \theta).
\end{align}
Here, the curve $x \colon \R \to \R^n$ is the trajectory of the plant starting from initial condition $x(0) =x_0$, and $N \in \N$ is the number of time steps in each rollout. Probing noise can be added to the input to encourage exploration, e.g, by instead setting $u_k = \hat{u}(x_k, \theta) + w_k$, where $w_k \sim \mathcal{N}(0, \sigma_w^2 I)$ is zero mean random noise.  Note that in the special case that $N=1$ the cost incurred when solving $\tilde{\textbf{P}}_\lambda$ approaches the cost incurred when solving $\textbf{P}_\lambda$ as $\Delta t \to 0$. In future work we plan to more formally study the relationship between these two problems. The policy optimization problem \eqref{eq:rl_prob} is in a standard form for reinforcement learning problems \cite{sutton2018reinforcement}, and in the following section we demonstrate how these discrete-time approximations can be used to successfully learn stabilizing controllers for unknown systems.


\section{Examples}
\label{sec:examples}
\subsection{Double Pendulum}
\begin{figure*}[h!]
\centering
\includegraphics[width=0.9\textwidth, trim=0.1cm 0cm 0.1cm 0cm, clip=true]{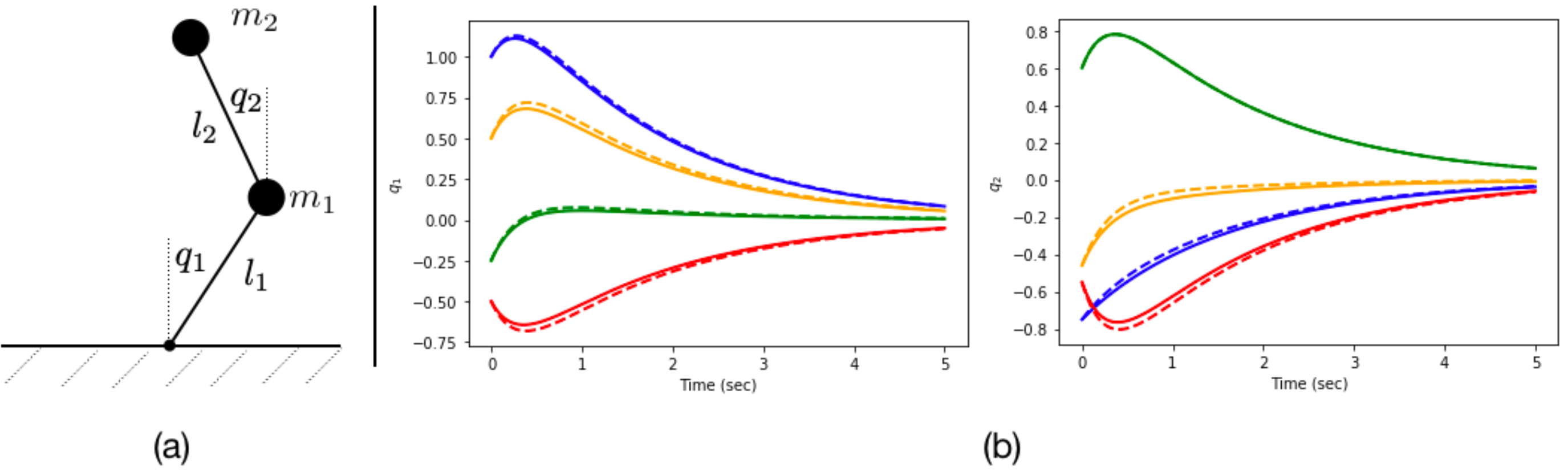}
\caption{(a) Depiction of the double pendulum model with the states and physical parameters shown. (b) Trajectories corresponding to different initial conditions for the learned controller and true min-norm controller for the system. Each color represents trajectories starting from a specific initial condition. Solid lines denote the trajectories generated by the true min-norm controller for the system while the dashed lines correspond to the trajectories generated by the learned controller. Observe that the learned controller closely matches the desired closed-loop behavior. Note that the velocities of the trajectories are not depicted, which is why several of the plotted curves intersect.}
\label{fig:dp}
\vspace{-12pt}
\end{figure*}

We first use our approach to learn a controller which stabilizes the double pendulum depicted in Figure~\ref{fig:dp} to the upright position, using the input-output linearization-based CLF design approach introduced in \cite{ames2014rapidly} to design the candidate CLF for the learning problem. The system has two generalized coordinates $q = (q_1,q_2)$ which represent the angles that each of the arms make with the vertical, with a motor attached at each joint so the system is fully actuated by torques $\tau = 
(\tau_1,\tau_2)$. The precise dynamics of the model can be found in \cite{capello2015modeling}, with states $x = (q,\dot{q})$ and input $u = \tau$. Below we use the configuration variables as outputs when applying the linearization-based CLF design \cite{ames2014rapidly}. 

As depicted in Figure \ref{fig:dp}a) the system is parameterized by the masses of the two arms $m_1, m_2$ as well as their lengths, $l_1, l_2 \in \R$. For the purposes of simulation, we set $m_1 = m_2 = l_1 = l_2 =1$. To set up the learning problem, we assume that we are given an inaccurate dynamics model with inaccurate estimates $\hat{m}_1, \hat{m}_2, \hat{l}_1, \hat{l}_2$. Specifically, we set $\hat{m}_1 = \hat{m}_2 = \hat{l}_1 = \hat{l}_2 =\frac{1}{2}$ so that each of the parameter estimates are half of their true value.

Using the input-output linearization design technique from \cite{ames2014rapidly}, we design a CLF for the system of the form
 $V \colon \R^4 \to \R$, $V := x^T P x$, with 
 \vspace{-5pt}
\begin{equation}
   P = \begin{bmatrix}
   1.5 I & 0.5I \\
   0.5I & 0.5I
   \end{bmatrix},
\end{equation}
where $I$ is the $2 \times 2$ identity matrix and by setting the desired dissipation rate to be $\sigma(x) = x^Tx$. This can be shown to be a valid CLF for both the inaccurate dynamics model and the true plant. We focus on learning the min-norm controller for the plant on the set $W^c = \set{ V(x) \leq c}$ with the design parameter $c =2$ and construct our learned controller by setting
\begin{equation}
    \hat{u}(x, \theta) = u_m(x) + \delta u(x,\theta) ,
\end{equation}
where $u_m$ is the min-norm CLF controller computed using the inaccurate dynamic parameters and the learned augmentation $\delta u$ is comprised of a linear combination of 500 radial basis functions so as to match the assumptions of Lemma \ref{lemma:convex}. 

\begin{figure}
\centering
\includegraphics[width=0.4\textwidth]{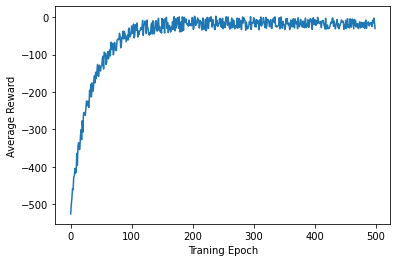}
\caption{Learning curve for the double pendulum.}\label{fig:dp_training}
\vspace{-10pt}
\end{figure}

We trained the learned component using a policy-gradient algorithm with action conditioned baselines \cite{sutton2018reinforcement}. Each training epoch consisted of $50$ 1-step roll-outs and a total of $500$ epoch were used. The time-step for the simulator was $0.05$ seconds. The performance of the ultimate learned controller is depicted in Figure~\ref{fig:dp}, where we see that the learned controller closely matches the behavior of the true min-norm controller for the system. To further evaluate the performance of the learned controller, we randomly selected 1000 states $\set{x_i}_{i=1}^{1000}$ in $W^c$ and calculated the ratio
\begin{equation}
    R = \sum_{i =1}^{1000} \frac{\norm{\hat{u}(x_i,\theta^*) - u_p^*(x_i)}_2 }{\norm{u_p^*(x_i)}_2},
\end{equation}
where $u_p^*$ is the true min-norm controller for the system and $\theta^*$ is the parameter selected by the training process. We calculated $R=0.044$, indicating that the learned controller was able to closely match the performance of the true min-norm controller for the system. As depicted in Figure~\ref{fig:dp_training}, the learning converges in about 200 iterations, which corresponds to about eight minutes of data. Our implementation of the learning algorithm for this problem was hand-coded, and we believe the sample efficiency for this problem could match that of the walking example below by improving the implementation.

\subsection{Bipedal Walking}

Next, we discuss how to apply our method to the Hybrid Zero Dynamics (HZD) framework using the CLF-based design approach proposed in \cite{ames2014rapidly} in order to learn an efficient, stable walking controller for a bipedal robot. We model the robot as a hybrid system with impulse effects as in \cite{ames2014rapidly},
\begin{equation} \label{eq:hybrid_biped_model}
  \Sigma:  
  \begin{cases}
      \dot{\eta} = f(\eta,z) + g(\eta,z)u, \\
      \dot{z} = h(\eta,z) &\text{when $(\eta,z) \notin \mathcal{S}$},\\
      \eta^+ = \Delta_X\left(\eta^-, z^-\right),\\
      z^+ = \Delta_Z\left(\eta^-, z^-\right)&\text{when $(\eta,z) \in \mathcal{S}$},
  \end{cases}
\end{equation}
 where $\eta \in \mathcal{X} \subset \R^{n_a}$ represents the controlled (actuated) states, $z \in \mathcal{Z} \subset \R^{n_u}$ represents the uncontrolled states and $u \in \mathcal{U} \subseteq \mathbb{R}^m$ represents the control inputs. The model assumes alternating phases of single support, where one foot is off the ground (swing foot) and the other (stance foot) is assumed to remain at a fixed point without slipping. The impact between the swing foot and the ground is modelled as a rigid impact and 
 occurs when $(\eta,z) \in \mathcal{S}$, where $\mathcal{S}$ is a smooth switching manifold. Here, $\eta^+ \in \mathcal{X}$ and $z^+ \in \mathcal{Z}$ represent the post-impact states while $\eta^- \in \mathcal{X}$ and $z^- \in \mathcal{Z}$ denote the pre-impact states. 

Following the framework in \cite{ames2014rapidly}, an input-output linearization based CLF is designed for the actuated coordinates during the continuous portion of the evolution of the state. Namely, we design a Lyapunov function $V \colon \R^{n_u} \to \R$ and dissipation rate $\sigma \colon \R^{n_u +n_a} \to \R_{\geq 0}$ such that the following condition holds for each $(\eta,z) \in \mathcal{X} \times \mathcal{Z}$: 
\begin{equation}
    \inf_{u \in U} \nabla V(\eta)[f(\eta,z) + g(\eta, z)u] \leq -\sigma(\eta, z).
\end{equation}
Thus, the control objective is to drive only the actuated states to zero. As shown in \cite{ames2014rapidly}, when the coordinates for the actuated and unactuated portions of the system are chosen correctly the condition $\eta \to 0$ corresponds to the robot converging to a periodic walking gate. The CLF and dissipation rate are designed so that the actuated coordinates are driven to zero fast enough to overcome shocks to the system introduced by the switching condition. We refer the readers to \cite{ames2014rapidly} for more details on this procedure.

To accommodate this new objective, our goal is to learn a control law $u \colon \mathcal{X} \times \mathcal{Z} \times \Theta \to \R^m$ such that 
\begin{equation}
   \underbrace{\nabla V(\eta)[f_p(\eta,z) + g_p(\eta, z)u(\eta,z, \theta)] +  \sigma(\eta, z)}_{\coloneqq  \hat{\Delta}(\eta,z, \theta)} \leq 0
\end{equation}
for each $(\eta, z) \in \mathcal{X}\times \mathcal{Z}$ for our choice of learned parameters $\theta \in \Theta$. Here, $f_p$ and $g_p$ are the terms in true dynamics of the plant, which may differ from the nominal dynamics in \eqref{eq:hybrid_biped_model}. To modify our approach to this new setting, for each $\lambda \in \R_{\geq0}$ we now define the loss
\begin{equation}
    \hat{L}_\lambda(\theta) = \mathbb{E}_{(\eta, z) \sim X} \left[\norm{u(\eta, z, \theta)}_2^2 + \lambda H(\hat{\Delta}(\eta, z, \theta))\right],
\end{equation}
where $X$ is now the uniform distribution over $\mathcal{X}\times \mathcal{Z}$. Despite the fact that the CLF is defined only over the lower dimensional state $\eta$, the theoretical results from section \ref{subsec:theory} naturally extend to this case. 
Moreover, the techniques from section \ref{subsec:rl} can be used to find local minimizers of $\hat{L}_\lambda$.

\begin{figure}
\centering
\includegraphics[width=0.5\textwidth]{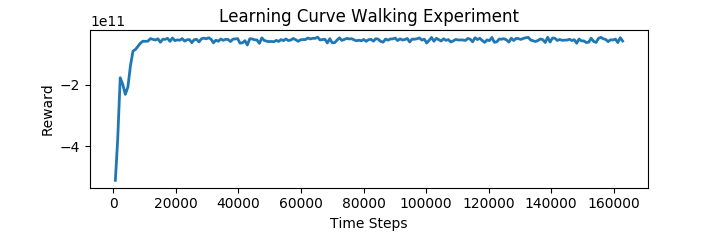}
\caption{Learning curve using PyBullet \cite{coumans2019} for the walking simulation.}\label{fig:walking-training}
\vspace{-10pt}
\end{figure}

 In particular, the proposed method is validated on a model for RABBIT \cite{chevallereau2003rabbit}, an under-actuated five-link planar bipedal robot with seven degrees-of-freedom. Model uncertainty is introduced by scaling the mass of each of RABBIT's links by a factor of two, i.e., the real plant's masses are twice the nominal model's masses. Our learned controlled is
\begin{align}
    \hat{u}(\eta,z,\theta) = u_m(\eta,z) + \delta u(\eta,z,\theta), 
\end{align}
where $u_m$ is the min-norm CLF controller obtained using the nominal model dynamics. The term $\delta u(\eta,z,\theta) \in \mathbb{R}^4$ takes the form of a Multi-Layer Perceptron (MLP) neural network with 2 hidden layers of width 64 each, tanh activation functions and layer normalization. 
We use the Soft Actor Critic algorithm \cite{levine_sac}, an off-policy method, for training the learned policy $\delta u(\eta,z,\theta)$. The training is done on episodes consisting of one walking step each. The simulations are conducted on the open-source physics simulator PyBullet \cite{coumans2019} using a discrete time-step of one millisecond. 
As it can be seen in Figure \ref{fig:walking-training}, the training converges in about 20,000 time steps, which corresponds to roughly 50 steps of the biped and about 20 seconds of data collection from the system. Altogether, the simulations and training took about 10 minutes of computation using the six cores of an Intel(R) Core(TM) i7-8705G CPU (3.10GHz), without using a GPU.

\begin{figure}
\centering
\includegraphics[width=0.5\textwidth]{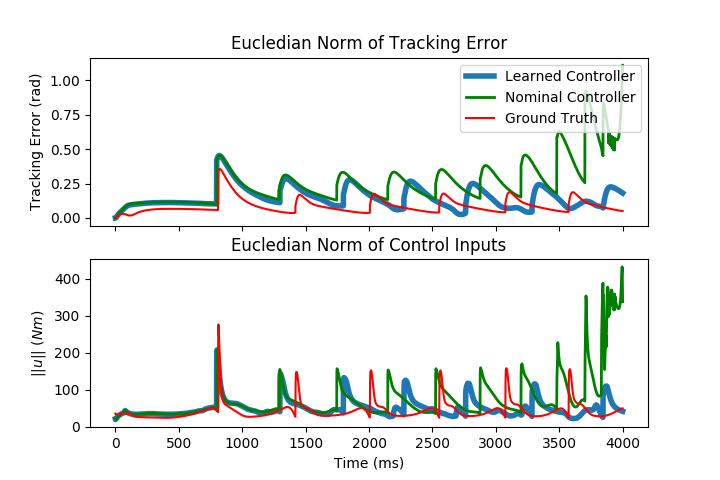}
\caption{Tracking error (top) and norm of control inputs (bottom) of the learned min-norm controller (blue), the nominal controller (green) and the actual CLF-based controller of the plant computed using the true robot dynamics (red), each simulated for 4 seconds of walking.}\label{fig:walking-compare-performance}
\vspace{-15pt}
\end{figure}

Figure \ref{fig:walking-compare-performance} shows a comparison between the proposed learned controller, the nominal controller $u_m$ and $u^*_p$, which is the CLF-based controller of the plant computed using the true (unknown) dynamics. This figure shows that while the nominal controller fails after ten walking steps making the robot fall, the learned controller 
achieves stable walking for an indefinite number of steps
and gives good tracking error performance. It is also important to notice that the learned controller achieves this while using similar magnitudes of control inputs as the nominal and the true CLF-based controllers. However, the tracking error performance is not as good as with the actual CLF-based controller of the plant, as expected, and is likely due to the fact that the learner has converged to a local minima. Additionally, we note that the walking speeds for the learned controller and the true min-norm CLF controller for the plant are different. Underactuated robots such as RABBIT may contain multiple periodic orbits on the surface $\{(\eta, z) \in \mathcal{X}\times \mathcal{Z} \colon \eta = 0 \}$. Thus, while both controllers successfully drive the system to this set, the periodic orbits the two controllers converge to are different.


\section{Conclusion}
\label{sec:conclusion}
There are several important avenues for future work. Demonstration of the proposed method on actual robotic systems in the face of sensor and actuator noise will be required before wider adoption of the technique. Future work will also investigate the inclusion of other local constraints on the dynamics of the closed loop system, such as those imposed by control barrier functions \cite{ames2019control}. On the theoretical side, sample complexity guarantees for the methods and a stronger characterization of the trade-off between our approach and model-based methods would more clearly characterize the advantages and disadvantages of the proposed approach. 

\appendix
\label{sec:proofs}
This Appendix contains proofs for several assertions made in the body of the document. 
\subsection{Proof of Lemma \ref{lemma:penalty}}

To prove the desired result, we demonstrate that for each $\theta^* \in \Theta \setminus \Xi$ there exists a finite $\bar{\lambda} \in \R_{\geq 0}$ such that $\theta^* \not \in S_{\lambda}$ for each $\lambda >\bar{\lambda}$. For a fixed $\theta^* \in \Theta\setminus \Xi$, define $M_1^{\theta^*} = E_{x\sim X}[\|\hat{u}(x,\theta^*)\|_2^2]$ and $M_2^{\theta^*} = E_{x \sim X}[H(\Delta(x,\theta^*))]$ so that for each $\lambda >0$ we have $L_\lambda(\theta^*) = M_1^{\theta^*} + \lambda M_2^{\theta^*}$. Since $\theta^* \not \in \Xi$, there must exist $x^* \in W^{c}$ such that $H(\Delta(x^*,\theta^*))>0$. Under our standing assumptions, the map $H(\Delta(\cdot, \theta^*))$ can be seen to be continuous, since the space of continuous functions is closed under addition, multiplication and composition. Putting these two facts together, there must exist a $\delta>0$ such that for each $x \in B^\delta(x^*) \cap W^{c}$ we have $H(\Delta(x, \theta^*))>0$. This in turn implies that $M_2^{\theta^*}>0$. Thus, we see that $L_\lambda(\theta^*) \to \infty$ as $\lambda \to \infty$.  

Next, letting $\bar{\theta}$ be defined as in the statement of the lemma, for each $\lambda \in \R_{\geq 0}$ we have $L_{\lambda}(\bar{\theta}) = M_1 ^{\bar{\theta}}$ where $M_1^{\bar{\theta}} = E_{x\sim X}[\|\hat{u}(x,\bar{\theta})\|_2^2]$ and we note that the term $E_{x \sim X}[H(\Delta(x,\theta^*))]$ contributes nothing to $L_{\lambda}(\bar{\theta})$ since $\bar{\theta} \in \Xi$. Thus, if we set $\bar{\lambda} = \max\set{0, \frac{M_1^{\bar{\theta}} - M_1^{\theta^*}}{M_2^{\theta^*}}}$ we see that $L_{\lambda}(\theta^*) >L_{\lambda}(\bar{\theta})$ for each $\lambda>\bar{\lambda}$, proving the desired statement for our fixed $\theta^*$. 
\subsection{Proof of Theorem \ref{thm:learn_min_norm}}
Let $\bar{\lambda}$ be defined as in the statement of Lemma \ref{lemma:penalty}. Then for each $\lambda > \bar{\lambda}$ we have $S_\lambda \subset \Xi$, where $\Xi$ is defined as in \eqref{eq:safe_values}. This implies that for each $\theta \in S_\lambda$ we have $L(\theta) = E_{x\sim X}[\|u(x,\theta)\|_2^2]$. Let $\bar{\theta}$ be defined as in the statement of the theorem, and let $\theta \in S_\lambda$ be arbitrary. By the definition of the min-norm control law we have $\|\hat{u}(x,\bar{\theta}) \|_2 \leq \|\hat{u}(x,\theta) \|_2$ for each $x \in W^{c}$, which in turn implies that $L(\bar{\theta}) \leq L(\theta)$. Next, suppose that $u(x^*, \theta) \neq u^*_p(x^*)$ for some $x^* \in W^{c}$. Again, using the definition of $u_p^*$ we have $\|\hat{u}(x^*,\bar{\theta}) \|_2 < \|\hat{u}(x^*,\theta) \|_2$. By the continuity of $\hat{u}(\cdot, \theta)$, we know that there exists $\delta > 0$ such that for each $x \in B^\delta(x^*)\cap W^c$ we have $\|\hat{u}(x,\bar{\theta}) \|_2^2 < \|\hat{u}(x,\theta) \|_2^2$. This implies that $L(\bar{\theta}) <L(\theta)$, demonstrating the desired result. 
\subsection{Proof of Lemma \ref{lemma:convex}}
To prove the claim, we will first consider the two maps $\theta \to E_{x \sim X} \left[\norm{u(x,\theta)}_2^2\right]$ and $\theta \to E_{x \sim X}\left[\lambda H(\Delta(x,\theta))\right]$ separately. In particular, we will show that the first term is strongly convex in $\theta$ while the second term is simply convex. The result of the theorem then follows from the fact that the addition of a strongly convex function and a convex function yields a strongly convex function. 

First, we rewrite $\norm{u(x,\theta)}_2^2$ as $\theta^TW(x)^TW(x)\theta$ where $W(x) = [u_1(x), u_2(x), \dots, u_{K}(x)]^T$ collects the basis of control functions. Note that the positive semi-definite matrix $\bar{W} = E_{x\sim X}\left[W(x)^TW(x)\right]$ is the Grammian for $\set{u_k}_{k=1}^K$ on $C(W^{c},\R^m)$, and thus will be full-rank and positive definite iff $\set{u_k}_k^K$ is linearly independent on this space. Based on these facts, we see that $E_{x \sim X} \left[\norm{u(x,\theta)}_2^2\right] = \theta^T \bar{W}\theta$ is a strongly convex quadratic function of the parameters. 

Next, we turn to the term $E_{x \sim X}\left[\lambda H(\Delta(x,\theta))\right]$. We demonstrate that for a fixed $x^* \in W^{c}$ and each $\lambda \in \R_{\geq 0}$ the mapping $\theta \to \|u(x^*,\theta)\|_2^2 + \lambda H(\Delta(x^*,\theta)) $ is strongly convex using basic properties of convex functions \cite{boyd2004convex}. We begin by examining the term $H(\Delta(x, \theta))$. Examining equations \eqref{eq:lin_param} and \eqref{eq:constraint_cost} we see that the map $\theta \to \Delta(x^*,\theta)$ is affine in $\theta$ for each fixed $x^* \in W^c$. Furthermore, we may rewrite the term $\lambda H(y) = \max \set{0, \lambda y}$. Since the pointwise maximum of two affine functions defines a convex function, we see that $\theta \to \lambda H(\Delta(x^*,\theta))$ is convex, implying that 
\vspace{-2pt}
\begin{equation*}
\lambda H(\Delta(x,\alpha \theta_3)) \leq \alpha \lambda H(\Delta(x, \theta_1)) + (1-\alpha) \lambda H(\Delta(x,
\theta_2))
\vspace{-3pt}
\end{equation*}
for each $x \in W^{c}$, $\theta_1,\theta_2 \in \R^K$ and $\theta_3 = \alpha \theta_1 + (1-\alpha)\theta_2$ for some $\alpha \in \sp{0,1}$. This pointwise fact implies that 
\vspace{-2pt}
\begin{align*}
E_{x \sim X} \left[\lambda H(\Delta(x,\theta_3))\right] &\leq \alpha E_{x \sim X} \left[\lambda H(\Delta(x, \theta_1))\right]\\
& +(1-\alpha) E_{x \sim X}\left[\lambda H(\Delta(x,
\theta_2))\right].
\end{align*}
Thus, $\theta \to E_{x \sim X}\left[\lambda H(\Delta(x,\theta))\right]$ is convex, as desired.

\bibliographystyle{IEEEtran}
\bibliography{IEEEabrv,refs.bib}

\end{document}